\UseRawInputEncoding
\documentclass{amsart}
\usepackage{amsmath}
\usepackage{amssymb}
\newtheorem{thm}{Theorem}[section]
\newtheorem{cor}[thm]{Corollary}
\newtheorem{conj}[thm]{Conjecture}
\newtheorem{lem}[thm]{Lemma}
\newtheorem{prop}[thm]{Proposition}
\newtheorem{cons}[thm]{Construction}
\theoremstyle{definition}
\newtheorem{defn}[thm]{Definition}
\theoremstyle{remark}
\newtheorem{rem}[thm]{Remark}
\newtheorem{ex}[thm]{Example}
\newtheorem{exs}[thm]{Examples}
\long\def\Thm#1{\begin{thm} #1 \end{thm}}
\long\def\Cor#1{\begin{cor} #1 \end{cor}}

\long\def\Lem#1{\begin{lem} #1 \end{lem}}
\long\def\Prop#1{\begin{prop} #1 \end{prop}}

\long\def\Rem#1{\begin{rem} #1 \end{rem}}
\long\def\Ex#1{\begin{ex} #1 \end{ex}}

\def\bar#1{\overline{#1}}
\def\Sect{\section}
\def\Rarr#1#2{\xrightarrow[#2]{#1}}
\def\Larr#1#2{\xleftarrow[#2]{#1}}
\def\Darr#1#2{{\scriptstyle #1}\downarrow{\scriptstyle #2}}

\long\def\Ref#1#2#3#4#5#6{
\bibitem{#1}
{\rm #2,}
\textit{#3.}
{\rm #4}
\textbf{#5}
{\rm #6.}
}
\long\def\Refb#1#2#3#4{
\bibitem{#1}
{\rm #2,}
\textit{#3.}
#4.
}
\def\Zz{{\mathbb Z}}
\def\Rr{{\mathbb R}}
\def\Cc{{\mathbb C}}
\def\Hh{{\mathbb H}}
\def\Ff{{\mathbb F}}
\def\Oo{{\mathbb O}}
\def\Kk{{\mathbb K}}

\def\phi{\varphi}
\def\into{\hookrightarrow}
\def\iso{\cong}
\def\leq{\leqslant}
\def\geq{\geqslant}

\def\comp{\mathbin{\mathchoice
{\circ}
{{\scriptstyle\circ}}
{{\scriptscriptstyle\circ}}
{{\scriptscriptstyle\circ}}
}}
\def\st{\mid}

\def\Hom{{\rm Hom}}
\def\cl#1{\bar{#1}}

\def\e{{\rm e}}
\def\i{{\rm i}}

\def\ss{\mathfrak{s}}
\begin{document}

\title[Fibrewise topological complexity of sphere and projective bundles]{On the fibrewise topological complexity of sphere and projective bundles}

\author{M.~C.~Crabb}
\address{%
Institute of Mathematics\\
University of Aberdeen \\
Aberdeen AB24 3UE \\
UK}

\email{m.crabb@abdn.ac.uk}
\begin{abstract}
We establish a stable homotopy-theoretic version of a recent result of
Farber and Weinberger \cite{farber1} on the fibrewise topological complexity of sphere bundles and prove, by closely parallel methods,
a similar result for 
real, complex and quaternionic projective
bundles. The symmetrized invariant introduced by Farber and Grant 
\cite{grant1} is also considered.
\end{abstract}
\maketitle
\section*{Introduction}
Two points $u,\, v\in S(V)$ in the unit sphere in a Euclidean space
$V$ (of dimension greater than $1$)
are joined by a unique shortest geodesic (in the standard
Riemannian metric) unless they are antipodal
and in that case the shortest geodesics
are parametrized by the unit sphere, of codimension $1$, 
in the orthogonal complement of the line $\Rr u=\Rr v$.

For a finite-dimensional $\Kk$-Hermitian vector space $V$ over
$\Kk =\Rr ,\, \Cc$ or $\Hh$, the points of the $\Kk$-projective
space $P_\Kk (V)$ are lines in $V$. Two lines $L,\ M\in P_\Kk (V)$
are joined by a unique shortest geodesic (in the metric determined
by the Hermitian structure) unless they are orthogonal and in that
case the shortest geodesics are parametrized by the sphere, of dimension $0$, $1$ or $3$, in the real vector space of 
$\Kk$-linear homomorphisms $L\to M$.

These two observations explain the relation between Farber's notion
\cite{farber-1, farber}
of the topological complexity of a sphere or projective space and
the existence of sections of an associated sphere bundle.
The same is true for the fibrewise\footnote{ 
On the question of terminology, I follow Ioan James, as in \cite{ioan},
in preferring the systematic use of the epithet `fibrewise' rather `parametrized' for topology over a base.}
topological complexity, \cite{farber0},
of the sphere bundle of a real vector
bundle or the projective bundle of a $\Kk$-vector bundle.
This relation is described for sphere bundles in
Theorem \ref{sphere}, following \cite{farber1},
and for projective bundles in Theorem \ref{proj}.
Symmetrized ($\Zz /2$-equivariant) 
versions, in the sense of \cite{grant1, grant2},
are given in Theorems \ref{symm} and \ref{symproj}.

The existence of a section of a sphere bundle is determined
in a stable range by the stable cohomotopy Euler class of
the vector bundle. Relevant concepts and notation are summarized
in Section \ref{Euler}.

Throughout the paper we shall use the notation
$t\mapsto c (t,u,v) : [0,1]\to S(V)$, with
$c (0,u,v)=u$, $c (1,u,v)=v$,
for the unique shortest geodesic joining two points $u,\, v\in S(V)$ in the unit sphere with $v\not=-u$.
To be precise,
$c (t,u,u)=u$,
and if $v\not=\pm u$ the $2$-dimensional
real vector space $\Rr u +\Rr v$ can be given a
complex structure such that $v=\e^{\i\theta}u$,
$0<\theta <\pi$ and then 
$c(t,u,v)=\e^{\i\theta t}u$.
We shall sometimes use without comment the fact that,
if $g : V\to V$ is an isometry, then 
$c(t,g\cdot u,g\cdot v)=g\cdot c(t,u,v)$.

\Sect{\label{Euler}
Preliminaries on the stable cohomotopy Euler class}
We shall use a notation for stable cohomotopy theory and
the stable cohomotopy Euler class that follows the classical
notation for cohomology and the cohomology Euler class.
For further details
we refer to \cite[II, Section 4]{ioan}.

Let $X$ be a compact
ENR\footnote{Euclidean Neighbourhood Retract.
In practice, $X$ usually admits
the structure of a finite complex and the sub-ENR 
$A$ is a sub-complex.}
and $A\subseteq X$ a closed sub-ENR. 
If $\alpha$ and $\beta$ are finite-dimensional real vector bundles
over $X$, we write $\omega^*(X,A;\, \alpha -\beta )$ for the
reduced stable cohomotopy group of the Thom space of the virtual
bundle $\alpha -\beta$.
There is a Hurewicz homomorphism to $\Zz$-cohomology, which we may
write using the Thom isomorphism as
$$
\omega^*(X.A;\, \alpha -\beta )\to H^{*-n}(X,A;\, \Zz (\alpha -\beta )),
$$
where $n=\dim\alpha -\dim\beta$ and $\Zz (\alpha -\beta )$ is the
local coefficient system of integers twisted by the orientation
bundle of $\alpha -\beta$.

Now consider an $n$-dimensional real vector bundle $\zeta$ over $X$
with sphere bundle $S(\zeta )$ and closed unit disc bundle $D(\zeta )$
(for a chosen Euclidean structure).
The Thom space of the pullback of $\zeta$ to $(D(\zeta ),S(\zeta ))$
is naturally identified with the Thom space of the virtual bundle
$\zeta -\zeta$ over $X$ and
there is a tautological Thom class
$$
u_\zeta\in \omega^0(D(\zeta ),S(\zeta );\, -\zeta),
$$
where, to be precise, the `coefficient system' $-\zeta$ is lifted from
$X$ to $D(\zeta )$.

The {\it stable cohomotopy Euler class} of $\zeta$
$$
\gamma (\zeta )\in\omega^0(X;\, -\zeta )
$$
is the restriction
of $u_\zeta$ to the zero-section 
$(X,\emptyset )\into (D(\zeta ),S(\zeta ))$, just as the cohomology
Euler class $e(\zeta )\in H^n(X;\, \Zz (-\zeta ))$ is the
restriction of the cohomology Thom class,
and $e(\zeta )$ is the Hurewicz image of $\gamma (\zeta )$.
If $s$ is a section of the restriction $S(\zeta\, |\, A)$
over $A$, the {\it relative stable cohomotopy Euler class}
$$
\gamma (\zeta ; s)\in\omega^0(X,A;\, -\zeta )
$$
is defined to be $\tilde s^*(u_\zeta )$, where
$\tilde s : X \to D(\zeta )$ is any extension of $s$
to a section of $D(\zeta )$ over $X$.
(Two such extensions $\tilde s$ are homotopic through a linear
homotopy.)
From the definition, it is clear that $\gamma (\zeta ; s)=0$
if $s$ extends to a section of the sphere bundle
$S(\zeta )$ over $X$. The converse is true in the
(meta) stable range $\dim X <2(n-1)$, essentially as a consequence
of Freudenthal's suspension theorem.

The Thom class $u_\zeta$ can itself be expressed as
 a relative Euler class, namely
$u_\zeta =\gamma (\zeta ; s)$ where $s$ is the diagonal section of the
pullback $S(\zeta )\times_B S(\zeta )$ through
$S(\zeta )\to X$ of $S(\zeta )$.
\Sect{Sphere bundles}
Let $V$ be a Euclidean vector space of dimension $n+1$.
Two points $u,\, v\in S(V)$ in the unit sphere with $v\not=-u$
are joined by a unique shortest geodesic, which we write
symmetrically as $t\in [-1,1]=D(\Rr )\mapsto \rho (t,u,v)$
with $\rho (1,u,v)=u$ and $\rho (-1,u,v)=v$:
$$
\rho (t,u,v)=c((1-t)/2, u,v).
$$

The shortest geodesics joining $u$ and $-u$ are
parametrized by the sphere $S((\Rr u)^\perp )$ in the orthogonal 
complement of the line $\Rr u$ in $V$. For $w\in S((\Rr u)^\perp)$
we write $\sigma (w;t,u,-u)$, $-1\leq t\leq 1$, for
the geodesic from $u$ to $-u$ passing through $w$:
$\sigma (w; t,u,-u)=c(t, w,u)$ for $0\leq t\leq 1$,
$c(-t,w,-u)$ for $-1\leq t\leq 0$.

We denote the real projective space of $V$ by $P(V)$ and the Hopf line
bundle over $P(V)$, with fibre at a point $L\in P(V)$ the 
$1$-dimensional subspace $L\subseteq V$, by $\eta$. 
Thus $\eta$ is a subbundle of the trivial bundle
with fibre $V$; its orthogonal complement,
of dimension $n$, is denoted by $\zeta$.
Let us write
$$
\tilde P(V)=\{ (u,v)\in S(V)\times S(V)\st v=-u\}\, .
$$
It projects as a double cover of $P(V)$, $(u,v)\mapsto [u]=[v]$,
and is identified by
projection to the first factor with $S(V)$.
The lift of $\zeta$ to $\tilde P(V)$ or $S(V)$ is denoted by
$\tilde\zeta$.
We begin with an elementary geometric lemma that sets up a diffeomorphism between the complement $B(\tilde\zeta )$ of
the sphere bundle $S(\tilde\zeta )$ in the unit disc bundle
$D(\tilde\zeta)$ and the complement of the diagonal
$\Delta (S(V))$ in $S(V)\times S(V)$.
\Lem{\label{diag}
The map
$$
\pi = (\pi_+,\pi_-): (D(\tilde\zeta ),S(\tilde\zeta ))\to 
(S(V)\times S(V),\Delta (S(V))), 
$$
$$
((u,v),w)\mapsto (\pi_+((u,v),w),\pi_-((u,v),w))\qquad
$$
$$
\qquad\textstyle
=
(\frac{1-\| w\|^2}{1+\| w\|^2}u+\frac{2}{1+\| w\|^2}w,
\frac{1-\| w\|^2}{1+\| w\|^2}v+\frac{2}{1+\| w\|^2}w),
$$
where $(u,v)\in \tilde P(V)$, 
$w\in (\Rr u)^\perp =(\Rr v)^\perp$,
identifies the open ball $B(\tilde\zeta )$ with
the complement of the diagonal $\{ (u,u)\in S(V)\times S(V)\}$
in $S(V)\times S(V)$.
}
\begin{proof}
This is clear if we write $w=te$, where $0\leq t\leq 1$ and
$e\in S((\Rr u)^\perp) =S((\Rr v)^\perp)$.
The map takes $((u,v),te)$ to
$$
(\cos (\theta)u+\sin (\theta)e,
\cos (\theta)v+\sin (\theta)e),
$$
where $\cos (\theta )=(1-t^2)/(1+t^2)$, $\sin (\theta )=2t/(1+t^2)$,
$0\leq\theta\leq \pi /2$.
\end{proof}
For $n\geq 1$, 
let $\xi$ be an $(n+1)$-dimensional Euclidean vector bundle
over a connected, compact ENR
$B$. 
The notation introduced for the vector space $V$
extends naturally
to the vector bundle $\xi$.
We write $S(\xi )\to B$ for the sphere bundle of $\xi$
and $P(\xi )\to B$ for the real projective bundle.
The Hopf line bundle $\eta$ over $P(\xi )$ is defined
as a subbundle of the pullback of $\xi$, and we write
$\zeta$ for its orthogonal complement. 
Let $\tilde\zeta$ over $S(\xi )=\tilde P(\xi )$ denote the pullback
of $\zeta$.

We consider the fibre product $S(\xi )\times_B S(\xi )$.
Given $(u,v)\in S(\xi_x )\times S(\xi_x )$ 
in the fibre over $x\in B$ with $v\not=-u$,
$\rho (t,u,v)$, $-1\leq t\leq 1$,
is a path from $u$ to $v$ in the sphere $S(\xi_x)$.
For $u\in S(\xi_x )$ and $w\in S(\tilde\zeta_x)$, 
$t\mapsto \sigma (w; t,u,-u)$ is a path in $S(\xi_x)$ from
$u$ to $-u$.
The construction in Lemma \ref{diag} gives a map
$\pi : (D(\tilde\zeta ),S(\tilde\zeta ))\to 
(S(\xi)\times_BS(\xi),\Delta (S(\xi ))$ over $B$.

Most of the content of the first main result is contained
in the work of Farber and Weinberger \cite{farber1, farber},
but formulated rather differently.
In the statement, maps and sections are understood to be continuous.
\Thm{\label{sphere}
{\rm (See \cite{farber1, farber}).}
Consider the following conditions on the real vector bundle $\xi$
involving an integer $k\geq 0$.
\par\noindent
{\rm (1)}
The vector bundle $k\tilde\zeta =\Rr^k\otimes\tilde\zeta$ 
over $S(\xi )=\tilde P(\xi )$ admits a nowhere zero section.
\par\noindent
{\rm (2)}
There is an open cover $V_1,\ldots ,V_k$ of
$S(\xi)$ such that for each $i=1,\ldots ,k$ there is
a fibrewise map 
$\psi_i:D(\Rr )\times V_i \to S(\xi )=\tilde P(\xi )$ satisfying
$\psi_i(1,u)=u$ and $\psi_i(-1,u)=-u$ for $u\in V_i$.
\par\noindent
{\rm (3)}
There is an open cover $U_0,\ldots ,U_k$ of
$S(\xi)\times_B S(\xi)$ such that for each $i=0,\ldots ,k$ there is
a fibrewise map $\phi_i:D(\Rr )\times U_i \to S(\xi )$ satisfying
$\phi_i(1,u,v)=u$ and $\phi_i(-1,u,v)=v$ for $(u,v)\in U_i$,
and $\phi_i(t,u,u)=u$
for all $t\in D(\Rr )$ and $(u,u)\in U_i$.
\par\noindent
{\rm (4)}
The stable cohomotopy Euler class 
$\gamma (\tilde\zeta )\in \omega^0(S(\xi );\, -\tilde\zeta )=
\omega^0(\tilde P(\xi );\, -\tilde\zeta )$ satisfies
$\gamma (\tilde\zeta )^k=0$.

Then the condition {\rm (1)} implies {\rm (2)}, 
condition {\rm (2)} implies
{\rm (3)} and {\rm (3)} implies {\rm (4)}.
If $\dim B <(2k-1)n-2$, then {\rm (4)} implies {\rm (1)}.
}
\begin{proof}
Of course, (1) implies (4).
The converse is true in the (meta) stable range $\dim S(\xi )=
\dim B +n <2(kn-1)$,
that is, $\dim B <(2k-1)n-2$.
\par\noindent (1)$\implies$(2).
Suppose that $(s_1,\ldots ,s_k)$ is a nowhere zero
section of $k\tilde\zeta$.
Take $V_i=\{ u\in S(\xi ) \st s_i(u)\not=0\}$.
Define
$\psi_i(t,u)$, in terms of $w=s_i(u)/\| s_i(u)\|$, to be
the geodesic path $\sigma (w; t, u,-u)$.
\par\noindent (2)$\implies$(3).
Take $U_0=\{ (u,v)\in S(\xi )\times_B S(\xi )\st v\not=-u\}$
and, using the diffeomorphism $\pi= (\pi_+,\pi_-) : B(\tilde\zeta )
\to \{ (u,v)\in S(\xi )\times_B S(\xi )\st u\not=v\}$
of Lemma \ref{diag}, take $U_i=\pi (B(\tilde\zeta \, |\, V_i))$
for $i=1,\ldots ,k$.

Define $\phi_0(t,u,v)$ to be $\rho (t,u,v)$, and,
for $i=1,\ldots ,k$, define
$\phi_i(t,u,v)$ where $(u,v)=\pi ((u_0,v_0),w)$ (so $v_0=-u_0$), with
$(u_0,v_0)\in V_i$,
to be
$$
\phi_i(t,u,v)=\begin{cases}
\pi_-((u_0,v_0), (-2t-1)w)
&\text{if $-1\leq t\leq -1/2$,}\\
\psi_i(2t,u_0)
&\text{if $-1/2\leq t\leq 1/2$,}\\
\pi_+((u_0,v_0), (2t-1)w)
&\text{if $1/2\leq t\leq 1$.}
\end{cases}
$$
\par\noindent (3)$\implies$(4).
We use Lemma \ref{diag}
to make, for any integer $j$, the identification
$$
\pi^*:
\omega^0(D(\tilde\zeta ),S(\tilde\zeta );\,-j\tilde\zeta )
\Rarr{\iso}{}
\omega^0(S(\xi )\times_B S(\xi ),S(\xi );\, -j\tilde\zeta).
$$

Our strategy is to show that the Thom class 
$$
u\in\omega^0(D(\tilde\zeta ),
S(\tilde\zeta);\, -\tilde\zeta )=\omega^0(\tilde P(\xi ))\cdot u
$$ 
vanishes on each of the $k+1$ 
open sets $\pi^{-1}(U_i)$ that cover $D(\tilde\zeta)$. 
It will then follow (see, for example, \cite[II: Lemma 3.14]{ioan})
that 
$$
u^{k+1}=
\gamma (\tilde\zeta )^k\cdot u\in 
\omega^0 (D(\tilde\zeta ),S(\tilde\zeta );
\,-(k+1)\tilde\zeta )
$$
is zero, and hence that
$\gamma (\tilde\zeta )^k\in
\omega^0 (\tilde P(\xi );\, -k\tilde\zeta )$
is zero. 
(The Thom class $u$ is the relative Euler class
$\gamma (\tilde\zeta ;s)$ of the inclusion $s : S(\tilde\zeta )
\into \tilde\zeta$ over $S(\tilde\zeta )$. The diagonal
inclusion $S(\tilde\zeta )\into \Rr^{k+1}\otimes\tilde\zeta$
is homotopic through nowhere zero sections to the inclusion
of the first factor.)

We identify the pullback of $\xi$ to $\tilde P(\xi )$ with
$\Rr\oplus\tilde\zeta$ by $(u,w)\mapsto (1,w)$,
$(v,w)\mapsto (-1,w)$ over $(u,v)\in \tilde P(\xi )$.
The Thom class $u$ corresponds to the relative Euler class
$$
\gamma (\tilde \zeta ;s)\in
\omega^0 (X,A;\,-(\Rr\oplus\tilde\zeta ))
=\omega^0 (D(\tilde\zeta ),S(\tilde\zeta );\, -\tilde\zeta )
$$
of the nowhere zero section $s$
over $A=(D(\Rr)\times S(\tilde\zeta ))
\cup (S(\Rr )\times D(\tilde\zeta ))$
of the pullback
of $\Rr\oplus\tilde\zeta$ to
$X=D(\Rr )\times D(\tilde\zeta )$
given by 
$$
s(t,w)
=
\begin{cases}
(0,w)&\text{if $t\in D(\Rr )$, $w\in S(\tilde\zeta)$;}\\
(\frac{1-\| w\|^2}{1+\| w\|^2}t,\frac{2}{1+\| w\|^2}w)&
\text{if $t\in S(\Rr )$, $w\in D(\tilde\zeta )$.}
\end{cases}
$$
Over $D(\Rr )\times \pi^{-1}(U_i)$, $s$ extends to the nowhere zero
section $s_i$ supplied by $\phi_i$ as
$s_i(t,w)=\phi_i(t,\pi ((u,v),w))\in 
S(\xi_x)=S(\Rr\oplus\tilde\zeta_{(u,v)} )$,
where $u\in S(\xi_x)$ and $v=-u$.
\end{proof}
\Rem{Direct proofs of the implications (1)$\implies$(3),
by an explicit construction of geodesic paths, and 
(2)$\implies$(4) in Theorem \ref{sphere} may be illuminating.
}
\begin{proof}
\par\noindent (1)$\implies$(3).
Suppose that $(s_1,\ldots ,s_k)$ is a nowhere zero
section of $k\tilde\zeta$.
Take $U_0=\{ (u,v)\in S(\xi )\times_B S(\xi )\st v\not=-u\}$
and $U_i=\{ (u,v)\in S(\xi)\times_B S(\xi )\st s_i(u)\not=0,\, 
\langle u,v\rangle <0\}$
for $i=1,\ldots ,k$.
Notice that $(u,-u)\in U_i$ if $s_i(u)\not=0$. 

Define $\phi_0(t,u,v)$ to be $\rho (t,u,v)$.
For $i=1,\ldots ,k$, define
$\phi_i(t,u,v)$, in terms of $w=s_i(u)/\| s_i(u)\|$, to be
$\sigma (w; t, u,v)$.
\par\noindent (2)$\implies$(4).
By interpreting $\psi_i$ as a homotopy
we shall show that the Euler class
$y=\gamma (\tilde\zeta )\in\omega^0(S(\xi );\, -\tilde\zeta )$ 
restricts to zero on $V_i$.
So the $k$-fold product $y^k=
\gamma (\tilde\zeta )^k\in \omega^0(S(\xi );\,-k\tilde\zeta )$
must be zero.

Now the Euler class $\gamma (\tilde \zeta)$ corresponds under
the Thom suspension isomorphism
$$
\omega^0(S(\xi );\, -\tilde\zeta )=
\omega^0(S(\xi );\, \Rr -\xi )\iso
\omega^0((D(\Rr ),S(\Rr ))\times S(\xi );\, -\xi )
$$
to the relative Thom class $\gamma (\xi ; s)$ of the section
$s$ of the pullback of $\xi$ to $S(\Rr )\times S(\xi )$
taking the value $u$ at $(1,u)$, $-u$ at $(-1,u)$.
The map $\psi_i$ extends $s$ to a nowhere zero section
over $D(\Rr )\times V_i$.  Hence $y$ restricts to zero on
$V_i$ as claimed.
\end{proof}
\Rem{(\cite[Corollary 17]{farber1}). 
If $\xi$ admits a complex structure, then $\tilde\zeta$ admits a
nowhere zero section.
Indeed, the complex structure provides 
in the fibre of $S(\tilde\zeta )$ over
$u\in S(\xi_x )$, $x\in B$, 
the vector $\i u$, and in
$S(\xi_x )$ the path
$t\mapsto \e^{\pi\i t}u$, $0\leq t\leq 1$, from $u$ to $-u$.

More generally, if there is an open cover $(W_i)_{i=1}^k$ of $B$
such that the restriction of $\xi$ to each open set $W_i$ admits
a complex structure, then $k\tilde\zeta$ admits
a nowhere zero section.
}
\Rem{If $\dim B < (k-1)n$, then condition (1) holds for 
purely dimensional reasons, and the other conditions
(2), (3), (4) follow.
}
\Prop{Condition {\rm (4)} in Theorem \ref{sphere} is implied by
the following weaker form of {\rm (3)}.
\par\noindent
{\rm ($3'$)}
There is an open cover $U_0,\ldots ,U_k$ of
$S(\xi)\times_B S(\xi)$ such that for each $i=0,\ldots ,k$ there is
a fibrewise map $\phi_i:D(\Rr )\times U_i \to S(\xi )$ satisfying
$\phi_i(1,u,v)=u$ and $\phi_i(-1,u,v)=v$ for $(u,v)\in U_i$.
}
The {\it parametrized topological complexity} of the bundle
$S(\xi )\to B$, in the sense of \cite{farber0, farber, farber1},
is the smallest integer $k$ for which condition ($3'$) holds.

The proof below is a reworking of the argument in \cite{farber1}.
\begin{proof}
\par\noindent ($3'$)$\implies$(4).
Consider, for any integer $j$,
the split short exact sequence
$$
0\to
\omega^0(D(\tilde\zeta ), S(\tilde\zeta);\, -j\tilde\zeta)
\iso
\omega^0(S(\xi)\times_B S(\xi ),S(\xi );\, -j\lambda^*\tilde\zeta )
$$
$$
\to
\omega^0(S(\xi)\times_B S(\xi );\, -j\lambda^*\tilde\zeta)\Rarr{\Delta^*}{\Larr{}{\lambda^*}}
 \omega^0(S(\xi );\, -j\tilde\zeta )\to 0,
$$
where $S(\xi)$ is included as the diagonal $\Delta : S(\xi )\to
S(\xi )\times_B S(\xi)$ split by the projection to the first factor
$\lambda : S(\xi)\times_B S(\xi )\to S(\xi )$
and the isomorphism is given by the map $\pi$ of Lemma \ref{diag}.

As we have already noted in the proof
of the implication (3)$\implies$(4), the $\omega^0(S(\xi ))$-module 
$\omega^0(D(\tilde\zeta ),S(\tilde\zeta );\, -\tilde\zeta )$
is free on the canonical Thom class $u$,
and $u^{k+1}=\gamma (\tilde\zeta )^k\cdot u\in
\omega^0(D(\tilde\zeta ),S(\tilde\zeta );\, -(k+1)\tilde\zeta )$.
Let $x$, satisfying $\Delta^*x=0$,
denote the image of $u$ in 
$\omega^0(S(\xi)\times_B S(\xi );\, -\lambda^*\tilde\zeta)$.

Now $\phi_i$ defines a homotopy 
$(u,v)\in U_i\mapsto (u,\phi_i(1-2t,u,v))$:
$U_i\to S(\xi)\times_B S(\xi)$, $0\leq t\leq 1$,
between the diagonal map $\Delta\comp \lambda :(u,v)\mapsto (u,u)$
and the inclusion $(u,v)\mapsto (u,v)$, and
this homotopy respects the projection $\lambda$
to the first factor.
Since $\Delta^*x=0$, we see that
$x$ restricts to $0$ on $U_i$.
Hence the $(k+1)$-fold product
$x^{k+1}=0\in 
\omega^0(S(\xi )\times_B S(\xi );\, -(k+1)\tilde\zeta )$.
It follows that $u^{k+1}=0$ and so that $\gamma (\tilde\zeta )^k=0$.
\end{proof}
By passing from stable cohomotopy to cohomology 
we can bound the topological complexity.
Consider the $\Zz$-cohomology Euler class $e(\tilde\zeta )
\in H^n(S(\xi );\,\tilde\Zz )$, where $\tilde\Zz$ denotes the
local coefficient system twisted by the orientation bundle
of $\xi$, as the Hurewicz image of the stable cohomotopy
Euler class 
$\gamma (\tilde\zeta )\in\omega^0(S(\xi );\, -\tilde\zeta )$. 
The power $e(\tilde\zeta )^k$ lies in 
$H^{kn}(S(\xi );\,\tilde\Zz^{\otimes k})$, 
where the $k$-fold tensor power is equal to $\Zz$ if
$k$ is even, $\tilde\Zz$ if $k$ is odd.
Let $\beta : H^i(B;\,\Ff_2)\to H^{i+1}(B;\,\Zz )$
and $\tilde\beta : H^i(B;\,\Ff_2)\to H^{i+1}(B;\,\tilde\Zz )$
denote the mod $2$ Bockstein homomorphisms.
The Pontryagin class $p_m(\xi )\in H^{4m}(B;\,\Zz )$
is equal to $(-1)^mc_{2m}(\Cc\otimes\xi )$
(the Chern class of the complexification).
\Prop{If condition {\rm (4)} in Theorem \ref{sphere}
holds, then $e(\tilde\zeta )^k=0$. This may be expressed in
terms of the cohomology of $B$ as follows.
\par\noindent{\rm (a)}
If $n$ is odd, 
$(\tilde\beta w_{n-1}(\xi ))^k\in H^{kn}(B;\,\tilde\Zz^{\otimes k} )$
is divisible by $e(\xi )\in H^{n+1}(B;\,\tilde\Zz )$.
\par\noindent{\rm (b)}
If $n=2m$ is even and $k=2l$ is even, then
$p_m(\xi )^l\in H^{kn}(B;\,\Zz )$ is divisible by
$e(\xi )\in H^{n+1}(B;\,\tilde\Zz )$,
and, hence, $2p_m(\xi )^l=0$.
\par\noindent{\rm (c)}
If $n=2m$ is even and $k=2l+1$ is odd, then
$2p_m(\xi )^l=0$ and
$\tilde\beta (w_n(\xi )q)\in H^{kn}(B;\,\Zz )$
is divisible by $e(\xi )$ for any class
$q\in H^{2ln-1}(B;\,\Ff_2)$ such that $\beta (q)=p_m(\xi )^l$.
}
\begin{proof}
We have Gysin sequences
$$
H^{i-n-1}(B;\,\Zz )\Rarr{e(\xi )}{} H^i(B;\,\tilde\Zz )
\to H^i(S(\xi );\,\tilde\Zz )\to
H^{i-n}(B;\,\Zz )\Rarr{e(\xi )}{} H^{i+1}(B;\,\tilde\Zz )
$$
and
$$
H^{i-n-1}(B;\,\tilde \Zz )\Rarr{e(\xi )}{} H^i(B;\,\Zz )
\to H^i(S(\xi );\,\Zz )\to
H^{i-n}(B;\,\tilde\Zz )\Rarr{e(\xi )}{} H^{i+1}(B;\,\Zz )\, .
$$

\par\noindent (a).
Since $\tilde\beta w_{n-1}(\tilde\zeta )=e(\tilde\zeta)$
and $w_{n-1}(\tilde\zeta )$ is the lift of $w_{n-1}(\xi )$,
the assertion follows from the exact sequence.
\par\noindent (b).
If $n=2m$ is even, then $e(\tilde\zeta )^2=p_m(\tilde\zeta )$,
which is the lift of $p_m(\xi )$.
So $e(\tilde\zeta )^{2l}$ is the lift of $p_m(\xi )^l$.
From the Gysin sequence, $p_m(\xi )^l$ is divisible by
the $2$-torsion class $e(\xi )$.
\par\noindent (c).
Because $n$ is even, the Euler class $e(\tilde\zeta )$
maps to $2\in H^0(B;\,\Zz )$ in the exact sequence, 
and hence $e(\tilde\zeta )^{2l+1}$ maps to $2p_m(\xi )^l$.
So $2p_m(\xi )^l=0$ and from the Bockstein exact
sequence $p_m(\xi )^l$ lifts to a class $q\in H^{2ln-1}(B;\,\Ff_2)$.
Thus $e(\tilde\zeta )^k=e(\tilde\zeta )\beta (q)
=\tilde\beta (w_n(\tilde \zeta )q)$.
But $w_n(\tilde\zeta )$ is the lift of $w_n(\xi )$.
\end{proof}

\Rem{(\cite[Theorem 14]{farber1}). 
If in case (c), $e(\xi )=0$, there is another criterion. Then
$\tilde\beta (w_{2m}(\xi ))=0$ and so $w_{2m}(\xi )$ lifts
to an integral class $x\in H^{2m}(B;\,\tilde\Zz )$.
The reduction (mod $2$) of $e(\tilde\zeta )$ is
also equal to $w_{2m}(\xi )$. So $e(\tilde\zeta )=x+2y$
for some class $y\in H^n(S(\xi );\,\tilde\Zz )$.
Hence, the Gysin sequence splits as
$$
H^*(S(\xi );\,\tilde\Zz ) = H^*(B;\, \tilde\Zz )1\oplus 
H^{*-n}(B;\,\Zz )y\, .
$$
We have $e(\tilde\zeta )^k= p_m(\xi )^l x + 2p_m(\xi )^ly$.
Thus, $e(\tilde\zeta )^k=0$ if and only if $p_m(\xi )^lx=0$
and $2p_m(\xi )^l=0$.

For example, if $\xi$ admits a stable complex structure 
we may take $x$ to be the $m$th Chern class.
}
\Ex{(\cite[Theorem 8]{farber-1}).
Let $B=*$. 
If $n$ is odd, then condition (1) of Theorem \ref{sphere}
holds if $k\geq 1$, because $\Rr^{n+1}$ admits a complex structure,
and if $k=0$ condition (4) fails, because $e(\tilde\zeta )^0=1$.
If $n$ is even, then (1) holds if $k\geq 2$, for dimensional reasons,
and (4) fails if $k=1$, because $e(\tilde\zeta )^1=2
\in H^n(S(\Rr^{n+1});\,\Zz )=\Zz$.
}
\Prop{{\rm (\cite[Theorem 2]{farber2})}.
If condition {\rm (4)} in Theorem \ref{sphere} holds, then 
$w_n(\xi )^k\in
H^{kn}(B;\,\Ff_2)$ is divisible by $w_{n+1}(\xi )$.
}
\begin{proof}
From consideration of the $\Ff_2$-Gysin sequence,
the $\Ff_2$-cohomology Euler class $e(\tilde\zeta )^k\in 
H^{kn}(S(\xi );\,\Ff_2)$ is zero if and only if $w_n(\xi )^k\in
H^{kn}(B;\,\Ff_2)$ is divisible by $w_{n+1}(\xi )$.
\end{proof}
We look next at the symmetry of the paths $t\mapsto\phi_i(t,u,v)$
appearing in Theorem \ref{sphere}.
The group $\Zz /2$ acts freely on $S(\xi )=\tilde P(\xi )$
by the antipodal involution $u\mapsto -u$; the orbit space
is $P(\xi )$. The group
acts on $S(\xi)\times_B S(\xi)$ by interchanging
the two factors; the fixed subspace is the diagonal $S(\xi )$.

Investigation of the $\Zz /2$-symmetry was begun in \cite{grant1}
and continued in \cite{grant2}.
\Thm{\label{symm}
For $k>1$, consider the following conditions.
\par\noindent
{\rm (0)}
There is a vector bundle homomorphism $B\times\Rr^k\to \xi$
over $B$ with rank $>1$ at each point of $B$.
\par\noindent
{\rm (1)}
The vector bundle $k\zeta =\Rr^k\otimes\zeta$ 
over $P(\xi )$ admits a nowhere zero section.
\par\noindent
{\rm (2)}
There is an open cover $V_1,\ldots ,V_k$ of
$S(\xi)$ by $\Zz /2$-invariant subsets 
such that for each $i$ there is
a fibrewise map 
$\psi_i:D(\Rr )\times V_i \to S(\xi )=\tilde P(\xi )$ satisfying
$\psi_i(1,u)=u$ and $\psi_i(-t,-u)=\psi (t,u)$ for $u\in V_i$,
$-1\leq t\leq 1$.
\par\noindent
{\rm (3)}
There is an open cover $U_0,\ldots ,U_k$ of
$S(\xi)\times_B S(\xi)$ by $\Zz /2$-invariant subsets 
such that for each $i$ there is
a fibrewise map $\phi_i:D(\Rr )\times U_i \to S(\xi )$ satisfying
$\phi_i(1,u,v)=u$ and $\phi_i(-1,u,v)=v$ for $(u,v)\in U_i$,
$\phi_i(-t,v,u)=\phi_i(t,u,v)$,
and $\phi_i(t,u,u)=u$ for all $(u,u)\in U_i$ and $t\in D(\Rr )$.
\par\noindent
{\rm (4)}
The stable cohomotopy Euler class 
$\gamma (\zeta )\in \omega^0(P(\xi );\, -\zeta )$ satisfies
$\gamma (\zeta )^k=0$.

Then the condition {\rm (0)} implies {\rm (1)}, 
condition {\rm (1)} implies {\rm (2)},  {\rm (2)} implies
{\rm (3)} and {\rm (3)} implies {\rm (4)}.

If $\dim B <(2k-1)n-2$, then {\rm (4)} implies {\rm (1)}.

If either
{\rm (a)} $n>1$ and $\dim B < (k-1)(n+1)$
or 
{\rm (b)} $n=1$ and $\dim B <2(k-1)-2$, 
then {\rm (4)} implies {\rm (0)}.
}
If $k=1$, none of the conditions holds: (0) trivially,
(4) because the Stiefel-Whitney class $w_n(\zeta )$ is non-zero.
\begin{proof}
\par\noindent (0)$\implies$(1).
A section of $\xi$ over $B$ determines by orthogonal projection
(after lifting to $S(\xi )$) a section of $\tilde\zeta$.
Thus, $k$ sections $r_1,\ldots ,r_k$ of $\xi$ determine $k$
sections $s_1,\ldots , s_k$ of $\tilde\zeta$, and the
section $s=(s_1,\ldots ,s_k)$ of $k\tilde\zeta$ will be
nowhere zero if, at each point $x\in B$, the vector subspace
of the fibre $\xi_x$ spanned by $r_1(x),\ldots ,r_k(x)$ has dimension
greater than $1$.
\par\noindent (1)$\implies$(2)$\implies$(3)$\implies$(4).
The verification proceeds {\it mutatis mutandis} as in the proof
of the corresponding implications in Theorem \ref{sphere},
using $\Zz /2$-equivariant stable homotopy in the third case.
\par\noindent (4)$\implies$(1) if $\dim B < (2k-1)n-2$.
This follows at once, 
because we are in the stable range $\dim P(\xi ) < 2(kn-1)$.
\par\noindent (4)$\implies$(0) if either (a) or (b).
The proof, which presupposes some familiarity with fibrewise stable 
homotopy theory, is presented in Section \ref{rank}.
\end{proof}
\Prop{\label{diagonal}
{\rm (See \cite[Theorem 5.2]{grant2}.)}
Condition {\rm (3)} in Theorem \ref{symm} is equivalent to the condition
\par\noindent{\rm ($3'$)}
There is an open cover $U'_0,\ldots ,U'_k$ of
$S(\xi)\times_B S(\xi)$ by $\Zz /2$-invariant subsets 
such that for each $i$ there is
a fibrewise map $\phi'_i:D(\Rr )\times U'_i \to S(\xi )$ satisfying
$\phi'_i(1,u,v)=u$ and $\phi'_i(-1,u,v)=v$ for $(u,v)\in U' _i$,
$\phi'_i(-t,v,u)=\phi'_i(t,u,v)$,
for all $t\in D(\Rr )$.
}
Adapting the terminology of \cite{gonzalez, grant2},
the {\it fibrewise symmetrized topological complexity} of the bundle
$S(\xi )\to B$
is the smallest integer $k$ for which condition ($3'$) holds.

The key property of the bundle $S(\xi )\to B$ that we shall use
is that it is fibrewise uniformly locally contractible
\cite[II: Definition 5.16]{ioan}.
\begin{proof}
We show that ($3'$) $\implies$ (3).
First of all, choose an open cover $U_0,\ldots ,U_k$ 
by $\Zz /2$-equivariant sets such that
$\cl{U_i}\subseteq U'_i$.

For a fixed $i$, we construct $\phi_i$ as follows.
Since $\phi'_i(t,u,u)=\phi'_i(-t,u,u)$ for $(u,u)\in
\cl{U_i}$ and $t\in [0,1]$, there is an open $\Zz /2$-subset
$\Omega$ of $\cl{U_i}$ containing all the diagonal points
$(u,u)$ and having the property that
$\phi'_i(t,u,v)\not=-\phi'_i(-t,u,v)$ for
all $(u,v)\in \Omega$ and $t\in [0,1]$.
Choose a continuous $\Zz /2$-invariant function $\tau : 
\cl{U_i}\to [0,1]$ such that $\tau (u,u)=1$ for all $(u,u)\in \cl{U_i}$
and $\tau (u,v)=0$ if $(u,v)\notin\Omega$.
We can now define
$$
\phi_i(t,u,v)=\begin{cases}
\phi'_i(t,u,v)&\text{if $|t|\geq \tau (u,v)$,}\\
\rho (s, \phi_i'(\tau (u,v),u,v), \phi_i'(-\tau (u,v),u,v))
&\text{if $t=s \tau (u,v)$, $s\in [-1,1]$.}
\end{cases}
$$
This function $\phi_i$ has the desired properties.
\end{proof}
A similar argument shows that the property (3) in 
Theorem \ref{sphere}
or Theorem \ref{symm} is a fibre homotopy invariant.
\Prop{Suppose that $\xi''$ is a vector bundle
over $B$ such that the sphere bundle $S(\xi'' )$ 
is fibre homotopy equivalent to $S(\xi )$
and that $\xi''$ satisfies the condition:

\smallskip
\par\noindent
There is an open cover $U''_0,\ldots ,U''_k$ of
$S(\xi'')\times_B S(\xi'')$ such that for each $i=0,\ldots ,k$ there is
a fibrewise map $\phi''_i:D(\Rr )\times U''_i \to S(\xi'' )$ satisfying
$\phi''_i(1,u,v)=u$ and $\phi''_i(-1,u,v)=v$ for $(u,v)\in U''_i$,
and, for all $t\in D(\Rr )$ and $(u,u)\in U''_i$,
$\phi''_i(t,u,u)=u$.

\smallskip
\par\noindent
Then $\xi$ satisfies condition {\rm (3)} of Theorem \ref{sphere}. 
}
\begin{proof}
Let $f: S(\xi )\to S(\xi'')$ and $g: S(\xi '')\to S(\xi )$ be inverse
fibre homotopy equivalences and $h_t: S(\xi )\to S(\xi )$,
$0\leq t\leq 1$, a fibre homotopy from the identity $h_0$
to $h_1=g\comp f$.

Put $U'_i=(f\times f)^{-1}U''_i$ and define
$\phi'_i : U_i'\to S(\xi )$ by
$$
\phi_i'(t,u,v)=
\begin{cases}
h_{2(1+t)}(v)&\text{if $-1\leq t\leq -1/2$;}\\
g(\phi''_i(2t,f(u),f(v)))&\text{if $|t|\leq 1/2$;}\\
h_{2(1-t)}(u)&\text{if $1/2\leq t\leq 1$.}
\end{cases}
$$
Notice that $\phi'_i(t,u,u)=\phi'_i(-t,u,u)$ if
$(u,u)\in U'_i$.
The construction in the proof of Proposition \ref{diagonal},
without the equivariance, produces the required maps
$\phi_i$.
\end{proof}
\Rem{\label{linear}
Here are two examples in which the condition (0) of Theorem \ref{symm}
holds.
\par\noindent{(i).\ }
If $\xi$ admits a $2$-dimensional trivial subbundle,
then (0) holds with $k=2$.
\par\noindent{(ii).\ }
Another example appears in \cite[Example 20]{farber1}; take $B$ to
be the complex projective space $P_\Cc (\Cc^{k-1})$
of dimension $k-2$ on $\Cc^{k-1}$, $k\geq 2$,
and $\xi$ to be the
$3$-dimensional direct sum of the complex Hopf line bundle and
the trivial real line bundle $B\times\Rr$. 
Let $r_1,\ldots ,r_{k-1}$ be the
sections of the Hopf bundle given by the
coordinate functions on $\Cc^{k-1}$ and let $r_k$ be the constant
section $1$ of the real line bundle.
}
\Prop{\label{sw}
If condition {\rm (4)} of Theorem \ref{symm} holds,
then the $\Ff_2$-Euler class $e(\zeta )\in
H^n(P(\xi );\,\Ff_2)$ satisfies $e(\zeta )^k=0$. In
terms of Stiefel-Whitney classes this says that
$$
(T^n+w_1(\xi )T^{n-1}+\ldots + w_n(\xi ))^k\in H^*(B;\, \Ff_2)[T]
$$
is divisible in the polynomial ring by 
$T^{n+1}+w_1(\xi )T^n+\ldots +w_{n+1}(\xi )$.
}
\begin{proof}
The cohomology ring of the projective bundle is
$$
H^*(P(\xi );\,\Ff_2)=H^*(B;\, \Ff_2)[T]/(T^{n+1}+w_1(\xi )T^n+
\ldots +w_{n+1}(\xi ))\, .
$$
The Euler class $e(\zeta )$, satisfying
$e(\zeta )e(\eta )=e(\xi )=w_{n+1}(\xi )$ is given by
$T^n+w_1(\xi )T^{n-1}+\ldots + w_n(\xi )$.
\end{proof}
For small $k$ the criterion in Proposition \ref{sw}
can be made explicit.
Write $t=e(\eta )$ and $x_i =t^i+w_1(\xi )t^{i-1}+\ldots +w_i(\xi )$,
$i=0.\ldots ,n$,
in $H^*(P(\xi );\,\Ff_2)$, so that $x_i=tx_{i-1}+w_i(\xi )$
(with $x_{-1}=0$) and $x_0,\ldots ,x_n$ is a basis 
of $H^*(P(\xi );\,\Ff_2)$
as a free module over $H^*(B;\,\Ff_2)$.

We have 
$$
e(\zeta )=x_n;\quad e(\zeta )^2= 
w_n(\xi )x_n+w_{n+1}(\xi )x_{n-1};
$$
$$
e(\zeta )^3 =
(w_n(\xi )^2+w_{n-1}(\xi )w_{n+1}(\xi ))x_n
+w_n(\xi )w_{n+1}(\xi )x_{n-1} + w_{n+1}(\xi )^2x_{n-2}.
$$
So, if $n\geq 0$, $e(\zeta )\not=0$; if $n\geq 1$, 
$e(\zeta )^2\not=0$ unless $w_{n+1}(\xi )=0$ and
$w_n(\xi )=0$; if $n\geq 2$,
$e(\zeta )^3\not=0$ unless 
$w_{n+1}(\xi )^2=0$, $w_n(\xi )w_{n+1}(\xi )=0$ and
$w_n(\xi )^2+w_{n-1}(\xi )w_{n+1}(\xi )=0$.

If $w_{n+1}(\xi )=0$, then $e(\zeta )^k =w_n(\xi )^{k-1}x_n$
is non-zero if and only if $w_n(\xi )^{k-1}\not=0$.
\Ex{(\cite[Theorem 6.1]{grant2}).
If $B=*$ and $n\geq 1$, 
the conditions of Theorem \ref{symm} hold if and only if
$k\geq 2$.
}
\begin{proof}
Since $e(\zeta )\not=0$, we certainly require $k\geq 2$.
If $k=2$ and $n\geq 1$, the condition (0) clearly holds.
\end{proof}
\Sect{Projective bundles}
Our discussion of projective bundles will run closely parallel
to the account of sphere bundles in the previous section
and we shall sometimes use the same notation (especially,
$\rho $ and $\sigma$, $\eta$ and $\zeta$),
but with a new meaning, for corresponding concepts.

We consider projective spaces over $\Kk =\Rr,\,\Cc$ and $\Hh$,
and write $d=\dim_\Rr \Kk$.
Let $V$ now be a finite-dimensional (left) $\Kk$-vector
space with a Hermitian inner product $\langle -,-\rangle :
V\times V\to \Kk$ (such that $\langle zu,v\rangle =z\langle u,v\rangle$
and $\langle u,zv\rangle =\langle u,v\rangle \bar z$ for
$u,v\in V$, $z\in\Kk$, 
and $\langle v,u\rangle=\bar{\langle u,v\rangle}$). This inner product determines a 
Euclidean structure on $V$ with $\| u\|^2=\langle u,u\rangle$.

The $\Kk$-projective space will be written as $P_\Kk (V)$;
its points are $1$-dimensional $\Kk$-subspaces $L\subseteq V$.
The Hopf line bundle $\eta$ over $P_\Kk (V)$ is the (left) $\Kk$-line
bundle over $P_\Kk (V)$ with fibre $L$ at a point $L\in P_\Kk (V)$.
Its orthogonal complement $\zeta$ over $P_\Kk (V)$ has fibre
$L^\perp\subseteq V$ at $L$.
For a unit vector $u\in S(V)$, we sometimes write
$[u]\in P_\Kk (V)$ for the line $\Kk u$.

If $n=1$, so that $V$ has dimension $2$, the projective space
$P_\Kk (V)$ is a sphere of dimension $d$, as we shall discuss
later in Remark \ref{dim1}.

For two lines $L, M$ which are not orthogonal,
the points $L,M\in P_\Kk (V)$ in the projective space
are joined by a unique shortest geodesic, which we write
as $\rho (t, L, M)$, $-1\leq t\leq 1$. To be precise,
$\rho (t,L,L)=L$,
and if $L\not=M$, where $L=\Kk u$
and $M=\Kk v$, with $\| u\| =1=\| v\|$ and with
$v$ chosen, given $u$, so that $\langle u,v\rangle\in\Kk$ 
is real and positive (which we can achieve, because $L$ and $M$ are not orthogonal), then
$\rho (t,L,M)=[c((1-t)/2,u,v)]$,
$\rho (1,L,M)=L$ and $\rho (-1,L,M)=M$.
(Changing $u$ to $zu$ with $z\in\Kk$, $|z|=1$,
changes $v$ to $zv$.)
If $L$ and $M$ are orthogonal,
the shortest geodesics joining $L$ and $M$
in $P_\Kk (V)$ are
parametrized by the sphere $S(\Hom_\Kk (L , M))$; we write them
as $\sigma (a; t,L,M)$, where $a: L\to M$ is a
$\Kk$-linear isometry, $-1\leq t\leq 1$, $\sigma (a; 1,L,M)=L$
and $\sigma (a; -1,L,M)=M$:
$$
\sigma (a; t, L, M)= [c((1-t)/2, u, a(u)],\text{\ \ where
$L=\Kk u$, $\| u\| =1$,}
$$
that is, 
$\sigma (a; t, L,M)=(\sin (\pi (t+1)/4)+\cos (\pi (t+1)/4)a)L$.
The Euclidean structure on the $d$-dimensional
real vector space $\Hom_\Kk (L,M)$ is specified by
requiring that $\| a(u)\| = \| a\| \cdot\| u\|$ for $u\in L$. 
We write $a^*\in\Hom_\Kk (M,L)$ for the Euclidean dual of $a$.
Let $\tilde Q(V)\subseteq P_\Kk (V)\times P_\Kk (V)$ 
denote the space of orthogonal
pairs $(L,M)$ and $\tilde\alpha$ the
$d$-dimensional Euclidean vector bundle over $\tilde Q(V)$ with
fibres $\Hom_\Kk (L,M)$.
\Lem{\label{symdiag}
The map $\pi =(\pi_+,\pi_-)\, :$
$$
(D(\tilde \alpha ),S(\tilde\alpha )) \to (P_\Kk (V)\times P_\Kk (V),
\Delta (P_\Kk (V))),
\quad (L,M,a)\mapsto  ((1+a)L, (1+a^*)M),
$$
where $L,\ M\in P_\Kk (V)$, $L$ is orthogonal to $M$, 
$a\in D(\Hom_\Kk (L,M))$, 
restricts to a diffeomorphism 
$$
B(\tilde\alpha )\to P_\Kk (V)\times P_\Kk (V)-\Delta (P_\Kk (V))
$$
from the open ball to the complement of the diagonal.
}
\begin{proof}
If $a\in S(\Hom_\Kk (L,M))$, $a^*=a^{-1}$, $M=aL$
and $(1+a^*)M=(1+a^{-1})aL=(a+1)L$. So $\pi$ does map
$S(\tilde\alpha )$ into $\Delta (P_\Kk (V))$.

Given $(L,M,a)\in D(\tilde\alpha )$, we can choose
$u\in L$ and $v\in M$ such that $\| u\|=1=\| v\|$
and $a(u)=tv$ with $t\in\Rr$ and $0\leq t\leq 1$.
Then $(1+a)L=\Kk x$ and $(1+a^*)M=\Kk y$,
where $x=(u+tv)/\sqrt{1+t^2}$, $y=(v+tu)/\sqrt{1+t^2}$
and $\langle x,y\rangle =2t/(1+t^2)$.
(We see, again, that if $t=1$, then $x=y$.)

In the opposite direction, given two distinct lines in
$P_\Kk (V)$, we may write them as $\Kk x$ and $\Kk y$
with $\| x\| =1=\| y\|$ and $\langle x,y\rangle =2t/(1+t^2)$
for $0\leq t<1$. Then $u=(x-ty)\sqrt{1+t^2}/(1-t^2)$
and $v=(y-tx)\sqrt{1+t^2}/(1-t^2)$ are orthogonal unit
vectors.
\end{proof}
Suppose that $\xi$ is an $(n+1)$-dimensional Hermitian
$\Kk$-vector bundle over a compact ENR $B$. 
We consider the fibre product $P_\Kk (\xi )\times_B P_\Kk (\xi )$.
Given $(L,M)\in P_\Kk (\xi_x )\times P_\Kk (\xi_x )$,
$x\in B$, with $L$ and $M$ not orthogonal,
$\rho (t,L,M)$, $-1\leq t\leq 1$,
is a path from $L$ to $M$ in $P_\Kk (\xi_x)$.
For orthogonal $L, M\in P_\Kk (\xi_x )$ and 
$a\in S(\Hom_\Kk (L,M))$, 
$t\mapsto \sigma (a; t,L,M)$ is a path in $P_\Kk (\xi_x)$ from
$L$ to $M$.
Write $\tilde Q(\xi )\subseteq P_\Kk (\xi )\times_B P_\Kk (\xi )$ for
the space of orthogonal pairs $(L,M)$ and let $\tilde\alpha$
be the (orthogonal) $d$-dimensional
real line bundle over $\tilde Q(\xi )$ with fibre
$\Hom_\Kk (L, M)$ at $(L,M)$.
Projection to the first factor $\tilde Q(\xi )\to P_\Kk (\xi )$ describes
$\tilde Q(\xi )$ as the projective bundle of $\zeta$
over $P_\Kk (\xi )$.
\Thm{\label{proj}
{\rm (See \cite{farber, landweber}).}
Consider the following conditions.
\par\noindent
{\rm (0)}
There is a $\Kk$-linear vector bundle monomorphism
$\zeta \into k\eta=\Rr^k\otimes\eta$ over $P_\Kk (\xi )$.
\par\noindent
{\rm (1)}
The real vector bundle $k\tilde\alpha =\Rr^k\otimes\tilde\alpha$ 
over $\tilde Q(\xi )$ admits a nowhere zero section.
\par\noindent
{\rm (2)}
There is an open cover $V_1,\ldots ,V_k$ of
$\tilde Q(\xi)$ such that for each $i$ there is
a fibrewise map $\psi_i:D(\Rr )\times V_i \to P_\Kk (\xi )$ satisfying
$\psi_i(1,L,M)=L$ and $\psi_i(-1,L,M)=M$ for $(L,M) \in V_i$.
\par\noindent
{\rm (3)}
There is an open cover $U_0,\ldots ,U_k$ of
$P_\Kk (\xi)\times_B P_\Kk (\xi)$ such that for each $i$ there is
a fibrewise map $\phi_i:D(\Rr )\times U_i \to P_\Kk (\xi )$ satisfying
$\phi_i(1,L,M)=L$ and $\phi_i(-1,L,M)=M$ for 
$(L,M)\in U_i$ and, 
for all $t\in D(\Rr )$, $\phi_i(t,L,L)=L$ for $(L,L)\in U_i$.
\par\noindent
{\rm (4)}
 The stable cohomotopy Euler class 
$\gamma (\tilde\alpha )\in \omega^0(\tilde Q(\xi );\, -\tilde\alpha )$ satisfies
$\gamma (\tilde\alpha )^k=0$.

Then the condition {\rm (0)} implies {\rm (1)},
{\rm (1)} implies
{\rm (2)}, {\rm (2)} implies {\rm (3)} and condition
{\rm (3)} implies {\rm (4)}.
If $\dim B <(2k-2n+1)d-2$, then {\rm (4)} implies {\rm (1)};
if $\dim B <(2k-3n+2)d-2$, then {\rm (4)} implies {\rm (0)}.

Moreover, the conditions {\rm (1)} and {\rm (2)} are equivalent. 
}
The proofs of the implications 
(1)$\implies$(2)$\implies$(3)
follow closely the corresponding deductions in Theorem \ref{sphere}.
\begin{proof}
\par\noindent (0)$\implies$(1).
Suppose given $r_i : \zeta \to \eta$, $i=1,\ldots ,k$,
such that, for $L\in P_\Kk (\xi )$, $w\in S(\zeta_{L})$,
the vectors
$(r_1)_L(w)$, $\ldots$, $(r_k)_L(w)$ in $L$ are not all zero.
Define $s_i(L,M)$ for $(L,M)\in \tilde Q(\xi )$ to be
the $\Kk$-linear map $L\to M$ dual to the restriction
$M\to L$ of $(r_i)_{L}: L^\perp \to L$. 
Then the section $(s_1,\ldots ,s_k)$ of $k\tilde\alpha$
is nowhere zero.
\par\noindent (1)$\implies$(2).
Suppose that $(s_1,\ldots ,s_k)$ is a nowhere zero
section of $k\tilde\alpha$.
Take $V_i=\{ (L,M)\in \tilde Q(\xi ) \st s_i(L,M)\not=0\}$,
and define
$\psi_i(t,L,M)$, in terms of $a=s_i(L,M)/\| s_i(L,M)\|$,
to be $\sigma (a; t, L, M)$.
\par\noindent (2)$\implies$(3).
Take $U_0=\{ (L,M)\in P_\Kk (\xi )\times_B P_\Kk (\xi )\st \,
\text{$L$ and $M$ are not orthogonal}\}$
and, using the diffeomorphism $\pi= (\pi_+,\pi_-) : B(\tilde\alpha )
\to \{ (L,M)\in P_\Kk (\xi )\times_B P_\Kk (\xi )\st L\not=M\}$
of Lemma \ref{symdiag}, take $U_i=\pi (B(\tilde\alpha \, |\, V_i))$
for $i=1,\ldots ,k$.

Define $\phi_0(t,L,M)$ to be $\rho (t,L,M)$, and,
for $i=1,\ldots ,k$, define
$\phi_i(t,L,M)$ where $(L,M)=\pi ((L_0, M_0]),a)$, with
$(L_0,M_0)\in V_i$,
to be
$$
\phi_i(t,(L,M))=\begin{cases}
\pi_-((L_0, M_0), (-2t-1)a)
&\text{if $-1\leq t\leq -1/2$,}\\
\psi_i(2t,L_0,M_0)
&\text{if $-1/2\leq t\leq 1/2$,}\\
\pi_+((L_0,M_0), (2t-1)a)
&\text{if $1/2\leq t\leq 1$.}
\end{cases}
$$
\par\noindent (3)$\implies$(4).
As in the proof of Theorem \ref{sphere}, we shall show that
the Thom class 
$u\in\omega^0(D(\tilde\alpha ), S(\tilde\alpha );
\, -\tilde\alpha )$ vanishes on each of the sets
$\pi^{-1}(U_i)$ covering $D(\tilde\alpha )$.
Thinking of $u$ as the relative Euler class $\gamma (\tilde\alpha ;s)$
of the diagonal section $s$ of the pullback of $S(\tilde\alpha )\subseteq\tilde\alpha$, we use homotopy-theoretic
parallel translation 
(or, simply, path lifting when $\Kk =\Rr$)
to construct sections $s_i$ extending $s$ to $\pi^{-1}(U_i)$.

Let $\eta_i$ over $U_i$ be the pullback of $\eta$ by
the map $(L,M)\in U_i\mapsto \phi_i(0, L, M)\in P_\Kk (\xi )$.
Choose an isometric $\Kk$-isomorphism between
$(\phi_i)^*\eta$ and $D(\Rr )\times \eta_i$ over
$D(\Rr )\times U_i$ extending the identity over $\{ 0\}\times U_i$.
This gives an orthogonal `parallel translation' isomorphism 
$A_i(L,M):
L=\eta_{L}\to\eta_{M}=M$, for $(L,M)\in U_i$. And $A_i(M,L)
=A_i(L,M)^*$.
Now for $(L,M,a)\in \pi^{-1}(U_i)$, define 
$s_i(L,M,a)\in \tilde\alpha_{(L,M)}$ to be
$$
e(M,a^*)^{-1}\comp A_i((1+a)L),(1+a^*)M)\comp e(L,a)
\in \Hom (L,M),
$$
where $e(L,a) : L\to (1+a)L$ and $e(M,a^*) : M\to (1+a^*)M$
are the isomorphisms given by $1+a$ and $1+a^*$.
If $a\in S(\Hom (L,M))$, then $a^*=a^{-1}$ and this
composition is equal to $a$.

\smallskip

\par\noindent (4)$\implies$(1) if $\dim B < (2k-2n+1)d-2$.
The stable range $\dim \tilde Q(\xi )< 2(dk-1)$
is $\dim B + dn+d(n-1)<2dk-2$, that is, $\dim B < (2k-2n+1)d-2$.
\par\noindent (4)$\implies$(0) if $\dim B <(2k-3n+2)d -2$.
We have already observed that $\tilde Q(\xi )$ is the
projective bundle $P_\Kk (\zeta )$ over the first 
factor $P_\Kk (\xi )$.
The vector bundle $\tilde\alpha$ is $\Hom (\eta_1,\eta_2)$, where
$\eta_1$ is the lift of $\eta$ over $P_\Kk (\xi )$ and $\eta_2$
is the Hopf line bundle of the projective bundle.
In the stable range $\dim  P_\Kk (\xi ) < 2(k-n)d+2(d-1)$,
that is, $\dim B <(2k-3n+2)d -2$, condition (0) is equivalent to
the vanishing of the stable cohomotopy Euler class
of $\Hom (\eta_2,\Rr^k\otimes\eta_1 )$ over $P_\Kk (\zeta )$,
which is isomorphic by duality to $k\tilde\alpha$.

\smallskip

\par\noindent (2)$\implies$(1).
For $(L,M)\in V_i$,
we again use homotopy-theoretic parallel translation
from $L$ to $M$ along the path $\psi_i(t,L,M)$.
(If the maps $\psi_i$ were fibrewise smooth we could use honest
parallel translation.)

Choose continuous functions $\mu_i : \tilde Q(\xi )\to [0,1]$
such that the closure $K_i$ of the support of $\mu_i$ is contained in
$V_i$ and the open sets $\mu_i^{-1}(0,1]$ cover $\tilde Q(\xi )$.
Let $\eta_i$ be the pullback of $\eta$ by the map $K_i\to P_\Kk (\xi )$:
$(L,M)\mapsto \psi_i(0,L,M)$.
Then the pullback of $\eta$ by the restriction of $\psi_i:
D(\Rr )\times K_i \to P_\Kk (\xi )$ is isomorphic to $D(\Rr )\times\eta_i$
by an isomorphism that is the identity over $\{ 0\}\times K_i$.
This isomorphism gives an isomorphism $A_i(L,M)$
by `parallel translation' from 
$\eta_{L}$ to $\eta_{M}$, that is, a non-zero element
of $\tilde\alpha_{L,M}$. Multiplying $A_i(L,M)$
by $\mu_i(L,M)$, we get a
section $s_i$ of $\tilde\alpha$ that is non-zero over
$\mu_i^{-1}(0,1]$. (The argument is simpler if
$\Kk=\Rr$. For then $S(\xi)\to  P(\xi )$ is a
double cover and we can use the unique path-lifting.)
\end{proof}
\Rem{In the stable range, condition (0) is equivalent to the
existence of a $\Kk$-monomorphism from the pullback of
$\xi$ into $(k+1)\eta$ over $P_\Kk (\xi )$.
If $\dim \tilde Q(\xi )< dk$, that is, $\dim B < (k-2n+1)d$, then
(1) holds. If $\dim B =(k-2n+1)d$, the question is answered
by cohomology.
}
\Ex{(\cite[Corollary 13]{farber}).
Let $\Kk =\Rr$, $B=*$, $\xi =V$, $n=1,3$ or $7$. Then
(0) holds with $k=n$.
}
\begin{proof}
Take $V=\Cc ,\, \Hh$ or $\Oo$ (the Cayley numbers) with the 
inner product $\langle u,v\rangle = {\rm Re}(u\bar v)$.
Then $(u,v)\mapsto {\rm Im}(u\bar v)$ induces
an embedding $\zeta\into I\otimes\eta$, where the space
$I=\{ w\in V\st \bar w=-w\}$ of imaginary numbers has dimension $k$.
\end{proof}
\Ex{(\cite[Section 7]{farber}).
Let $\Kk=\Rr$, $B=*$, $\xi =V$, $n=2$. Then (0) holds with $k=3$.
}
\begin{proof}
Choosing an orientation of $V$ (of dimension $3$), 
we can use the vector product
$\times : V\times V\to V$ to write down
a vector bundle inclusion $\zeta\into V\otimes\eta$
over the real projective plane $P(V)$.
\end{proof}
The real vector bundle $\tilde\alpha$ over 
$\tilde Q(\xi )\subseteq P_\Kk (\xi )\times_B P_\Kk (\xi )$
is the restriction of a vector bundle $\hat\alpha$ over 
$P_\Kk (\xi)\times_B P_\Kk (\xi )$ with fibre $\Hom_\Kk (L,M)$ over
$(L,M)$. On the diagonal $\Delta (P_\Kk (\xi ))$, $\hat\alpha$ has an 
obvious nowhere zero section given by the identity
$1\in \Hom_\Kk (L,L)$ over $(L,L)$.
\Prop{Consider the following conditions on the vector bundle $\xi$ in
Theorem \ref{proj}.
\par\noindent
{\rm ($3'$)}
There is an open cover $U_0,\ldots ,U_k$ of
$P_\Kk (\xi)\times_B P_\Kk (\xi)$ such that for each $i$ there is
a fibrewise map $\phi_i:D(\Rr )\times U_i \to P_\Kk (\xi )$ satisfying
$\phi_i(1,L,M)=L$ and $\phi_i(-1,L,M)=M$ for 
$(L,M)\in U_i$.
\par\noindent
{\rm ($4'$)}
The stable cohomotopy Euler class 
$\gamma (\hat\alpha )\in \omega^0(P_\Kk (\xi )\times_B P_\Kk (\xi );\, -\hat\alpha )$ satisfies
$\gamma (\hat\alpha )^{k+1}=0$.

Then condition {\rm ($3'$)} implies {\rm ($4'$)}.

The condition {\rm ($4'$)} implies that the class
$\gamma (\tilde\alpha )\in\omega^0(\tilde Q(\xi );\, -\tilde\alpha )$
satisfies $\gamma (\tilde\alpha )^{k+1}=0$
and, if $\dim B <(k-n+1)d-1$, the condition
{\rm (4)} that $\gamma (\tilde\alpha )^k=0$.
}
\begin{proof}
($3'$)$\implies$($4'$).
Choose a partition of unity $(\mu_i)$ subordinate to the
cover $(U_i)$.
As in the proof that (3) imples (4) in Theorem \ref{proj},
we get parallel translation isomorphisms
$A_i(L,M) : L\to M$ for $(L,M)\in U_i$, that is,
a non-zero element of $\hat\alpha_{(L,M)}$.
Multiplying by $\mu_i$, we get a section $s_i$ of $\hat\alpha$
that is non-zero where $\mu_i$ is non-zero.
Then $(s_i)$ is a nowhere zero section of
$\Rr^{k+1}\otimes\hat\alpha$ and the Euler class must be zero.

If $\gamma (\hat\alpha )^{k+1}=0$, then evidently
$\gamma (\tilde\alpha )^{k+1}=0$, because
$\hat\alpha$ restricts to $\tilde\alpha$ on
$\tilde Q(\xi )\subseteq P_\Kk (\xi )\times_B P_\Kk (\xi)$.
The final assertion follows from
consideration of the exact sequence of the pair 
$(P_\Kk (\xi )\times_B P_\Kk (\xi ),\Delta (P_\Kk (\xi )))$:
$$
\omega^{-1}(P_\Kk (\xi );\, -(k+1)\hat\alpha )\to
\omega^0(P_\Kk (\xi )\times_B P_\Kk (\xi );\, -(k+1)\hat\alpha)
\to
\omega^0(P_\Kk (\xi );\, -(k+1)\hat\alpha  ),
$$
because the group $\omega^{-1}(P_\Kk (\xi );\, -(k+1)\hat \alpha )$
is zero if $\dim B +dn +1 < (k+1)d$.
\end{proof}
Let us write $R=\Ff_2$ if $\Kk =\Rr$, $R=\Zz$ if $\Kk =\Cc$ or
$\Hh$, and let $w_i^\Kk (\xi )\in H^{di}(B;\, R)$ denote
the Stiefel-Whitney class $w_i(\xi )$ if $\Kk=\Rr$,
the Chern class $c_i(\xi )$ if $\Kk=\Cc$
and $c_{2i}(\xi )$ if $\Kk =\Hh$.
Thus, if $\xi$ is a $\Kk$-line bundle,
$w_1^\Kk (\xi )\in H^d(B;\, R)$ is the cohomology Euler
class $e(\xi )$.
\Prop{The cohomology of $\tilde Q(\xi )$ with $R$-coefficients
is described as
$$
H^*(\tilde Q(\xi );\, R)=
$$
$$
H^*(P_\Kk (\xi );\, R)[T]/
(w_n^\Kk (\xi )+\sum_{i=1}^n
(-1)^i(T^{i} +\ldots +S^jT^{i-j}+ \ldots +S^i)w_{n-i}^\Kk (\xi))
$$
where
$$
H^*(P_\Kk (\xi );\, R)=H^*(B;\, R)[S]/
(w_{n+1}^\Kk (\xi )-Sw_n^\Kk (\xi )+\ldots +(-1)^{n+1}S^{n+1}).
$$
The Euler classes of the Hopf line bundles on the
two factors are given by $S$ and $T$,
and the Euler class $e(\tilde\alpha )$ by $T-S$.
}
\begin{proof}
As the projective bundle $P_\Kk (\zeta )\to P_\Kk (\xi )$
of the $n$-dimensional bundle $\zeta$,
the space $\tilde Q(\xi )$ has cohomology ring
$$
H^*(\tilde Q(\xi );\, R)=H^*(P_\Kk (\xi );\, R)[T]/
(w_n^\Kk (\zeta )-Tw_{n-1}^\Kk (\zeta )+ \ldots +(-1)^nT^n),
$$
We have $(1+S)(1+w_1^\Kk (\zeta )+\ldots + w_n^\Kk (\zeta ))=
1+w_1^\Kk (\xi )+\ldots +w_{n+1}^\Kk (\xi )$. 
So $w_i(\zeta )^\Kk =w_i^\Kk (\xi )-Sw_{i-1}^\Kk (\xi )
+\ldots +(-1)^iS^i$.
\end{proof}
Notice that $(T-S)(T^i+ST^{i-1}+\ldots +S^i)=T^{i+1}-S^{i+1}$.
This checks the symmetry in $S$ and $T$.
\Ex{\label{milnor}
(\cite[Theorem 6]{farber}, but
going back to 1957 lectures of Milnor \cite[Theorem 4.8]{milnor}.)
For $\Kk =\Rr$, $B=*$, and $n=2^r$, $r\geq 1$,
condition (0) in Theorem \ref{proj}
holds if $k\geq 2^{r+1}-1$,
but condition (4) fails if $k=2^{r+1}-2$.
}
\begin{proof}
Condition (0) holds for dimensional reasons if $k>2n-1$.
The cohomology ring is 
$\Ff_2[S,T]/(S^{n+1}, \, T^n+ST^{n-1}+\ldots +S^n)$.
For $k=2n-1$, $(T+S)^{2^{r+1}-1}=0$.
For $k=2n-2$,
$(T+S)^{2^{r+1}-2}=(T^2+S^2)^{2^r-1}
=T^{2^r}S^{2^r-2}+T^{2^r-2}S^{2^r}\not=0$.
\end{proof}
\Ex{(\cite[Corollary 2]{farber}).
For $\Kk =\Cc$ and $B=*$, condition (0)
in Theorem \ref{proj} holds if
$k\geq 2n$, but (4) fails if $k=2n-1$.
}
\begin{proof}
If $k >2n-1$, then (0) holds for dimensional reasons.
The cohomology ring is 
$\Zz [S,T]/(S^{n+1}, \, T^n+ST^{n-1}+\ldots +S^n)$.
For $k=2n-1$,
$(T-S)^{2n-1}=(-1)^n\binom{2n-1}{n}(T^{n-1}S^n-T^nS^{n-1})$
is non-zero.
\end{proof}
\Ex{
In the same way,
for $\Kk =\Hh$ and $B=*$, condition (0) in Theorem \ref{proj}
holds if $k\geq 4n$, but (4) fails if $k=4n-1$.
}
As for sphere bundles, we can look for a symmetric version
of the main theorem.
The group $\Zz /2$ acts on $\tilde Q(\xi )$
by interchanging the two factors; the orbit space
is $Q(\xi )$.
There is a compatible involution on $\tilde\alpha$
given by $* :\Hom_\Kk (L,M)\to\Hom_\Kk (M,L)$;
the quotient vector bundle over $Q(\xi )$ is denoted by
$\alpha$.
\Thm{\label{symproj}
Consider the following conditions.
\par\noindent
{\rm (1)}
The real vector bundle $k\alpha =\Rr^k\otimes\alpha$ 
over $Q(\xi )$ admits a nowhere zero section.
\par\noindent
{\rm (2)}
There is an open cover $V_1,\ldots ,V_k$ of
$\tilde Q(\xi)$ by $\Zz /2$-invariant subsets
such that for each $i$ there is
a fibrewise map $\psi_i:D(\Rr )\times V_i \to P_\Kk (\xi )$ satisfying
$\psi_i(1,L,M)=L$ and $\psi_i(-1,L,M)=M$ for 
$(L,M)\in V_i$ and $\psi_i(-t,L,M)=\psi_i(t,M,L)$
for all $t\in D(\Rr )$.
\par\noindent
{\rm (3)}
There is an open cover $U_0,\ldots ,U_k$ of
$P_\Kk (\xi)\times_B P_\Kk (\xi)$ by $\Zz /2$-invariant subsets
such that for each $i$ there is
a fibrewise map $\phi_i:D(\Rr )\times U_i \to P_\Kk (\xi )$ satisfying
$\phi_i(1,L,M)=L$ and $\phi_i(-1,L,M)=M$ for 
$(L,M)\in U_i$, and, for all $t\in D(\Rr )$,
$\phi_i(-t,L,M)=\phi_i(t,M,L)$ for $(L,M)\in U_i$
and $\phi_i(t,L,L)=L$ for $(L,L)\in U_i$.
\par\noindent
{\rm (4)}
 The stable cohomotopy Euler class 
$\gamma (\alpha )\in \omega^0( Q(\xi );\, -\alpha )$ satisfies
$\gamma (\alpha )^k=0$.

Then the condition 
{\rm (1)} implies
{\rm (2)}, {\rm (2)} implies {\rm (3)} and condition
{\rm (3)} implies {\rm (4)}.
If $\dim B <2d(k-n)+d-2$, then {\rm (4)} implies {\rm (1)}.

Moreover, the conditions {\rm (1)} and {\rm (2)} are equivalent. 
}

\begin{proof}
The implications 
(1)$\implies$(2) 
$\implies$(3)$\implies$(4), (4)$\implies$(1) in the stable range,
and (2)$\implies$(1)
follow very closely the proofs of the corresponding assertions
in Theorem \ref{proj}. In the deduction of (4) from (3)
we need to use $\Zz /2$-equivariant stable homotopy.
And in the deduction of (1) from (2) 
we must choose $\mu_i$, $K_i$ and
the isomorphism between the pullback of $\eta$
to $D(\Rr )\times K_i$ and $D(\Rr )\times\eta_i$ to be symmetric.
\end{proof}
\Rem{Condition (3) in Theorem \ref{symproj} is equivalent to
\par\noindent
{\rm (1$'$)}
{\it The pullback of $\Rr^{k+1}\otimes \alpha$ to $D(\alpha )$ admits
a nowhere zero section extending the diagonal inclusion $S(\alpha )\into S(\Rr^{k+1}\otimes\alpha )$ over $S(\alpha)$,}
\par\noindent
from which (4) readily follows.
}
\begin{proof}
\par\noindent
(3)$\implies$(1$'$). 
We resume the discussion in the proof of the implication 
(3)$\implies$(4) above as set out in the corresponding step
in the proof of Theorem \ref{proj}.
Choose a $\Zz /2$-equivariant partition of unity $\mu_i$ subordinate
to the cover  $(\pi^{-1}U_i)$ of $D(\tilde\alpha )$
and form the global section $(\mu_i s_i)$ of 
$\Rr^{k+1}\otimes\tilde\alpha$.
At a point of $S(\tilde\alpha )$ the value of $\mu_is_i$
is a non-negative multiple of the value of $s$.
So the restriction of $(\mu_i s_i)$ to $S(\tilde\alpha )$ is
linearly homotopic through nowhere zero sections to
the diagonal inclusion and can be deformed to a section
that coincides with the diagonal section on $S(\tilde\alpha )$.
\par\noindent
(1$'$)$\implies$(3). 
Conversely, suppose that we have such a section given by
an equivariant section $(s_i)$ of $(k+1)\tilde\alpha$. 
Let $U_i = \pi (\{ (L_0,M_0,a)\in D(\tilde\alpha )\st 
s_i(L_0,M_0,a)\not=0\}$. It is an open neighbourhood of
the diagonal $\Delta (P_\Kk (\xi ))\subseteq P_\Kk (\xi )\times_B
P_\Kk (\xi )$.
For $(L,M)\in U_i$, writing $(L,M)=((1+a)L_0,(1+a^*)M_0)$
where $s_i(L_0,M_0,a) : L_0\to M_0$ is non-zero, so an isomorphism,
define 
$$
\phi_i(t,L,M) = [c((1-t)/2, (1+a)u, (1+a^*)eu)],
$$
in terms of $e=\| s_i(L_0,M_0,a)\|^{-1} s_i(L,M,a)$
and a generator $u$ of $L_0=\Kk u$ with $\| u\|^2=(1+\| a\|^2)^{-1}$.
To see that this makes sense as a definition,
notice first that $(1+a)u\not=-(1+a^*)eu$. For if 
$au=-eu$, we have $a=-e$ and so $\| a\|=1$, and then $a=e$
because $s_i$ is the identity on $S(\alpha )$, which forces
$u\not =-a^*eu$.
Secondly, on the diagonal, where $L=M$,
the expression for $\phi_i$ gives $L$, independently of the
choice of $(L_0,M_0,a)$.
For then $a=e$, and so $(1+a)u=(1+a^*)eu$.
\end{proof}
For a $2$-dimensional $\Kk$-vector space $V$, let us write
$\ss (V)$ for the $(d+1)$-dimensional real vector space of
$\Kk$-Hermitian endomorphisms of $V$ with real trace zero.
If $V=\Kk\oplus\Kk$, elements of $\ss (V)$ can be written as
matrices 
$$
\left[\begin{matrix} b&\,\,  a\\ \bar a&-b\end{matrix}
\right], \quad
\text{where $b\in\Rr$, $a\in\Kk$.}
$$
This allows us to identify $\tilde Q(V)$ with the sphere
$S(\ss (V))$ by mapping $(L,M)$ to the endomorphism
$(1,-1): V=L\oplus M\to L\oplus M=V$,
and thus to identify 
$Q(V)$ with the real projective space $P(\ss (V))$.
Furthermore, given $(L,M)$, we have an isomorphism
$\Rr\oplus \Hom_\Kk (L,M) \to \ss (V)$:
$$
(b,a)\mapsto 
\left[\begin{matrix} b&\,\,  a\\  a^*&-b\end{matrix}
\right] :
L\oplus M\to L\oplus M.
$$
This identifies the $d$-dimensional real vector bundle 
$\alpha$ over $Q(V)$ with the orthogonal complement of
the (real) Hopf line bundle over $P(\ss (V))$ in the trivial bundle
$\ss (V)$.

\Rem{\label{dim1}
For a $2$-dimensional vector bundle $\xi$, when $n=1$, these
constructions set up a precise correspondence between
the projective bundle $P_\Kk (\xi )$ in this Section
and the sphere bundle
$S(\ss (\xi ))$ of Section 2, 
the bundle $\tilde Q(\xi)$ and the sphere bundle
$S(\ss (\xi ) )=\tilde P(\ss(\xi ))$, the bundle $Q(\xi )$ and
the real projective bundle $P(\ss (\xi ))$,
in which $\alpha$ corresponds to the orthogonal complement,
called $\zeta$ in Section 2,
of the real Hopf bundle over $P(\ss (\xi ))$.
}
Resuming the discussion of a bundle $\xi$ of arbitrary dimension $n+1$,
let $\lambda$ be the real line bundle over $Q(\xi )$ associated
with the involution.
There is a projection $Q(\xi )\to G_2^\Kk (\xi )$ to the Grassmann bundle
of $2$-dimensional $\Kk$-subspaces of the fibres of $\xi$ taking
a pair $\{ L,M\}$ to the $2$-dimensional subspace
$L\oplus M$.
And $Q(\xi )=P(\ss (\beta))$ is the real projective bundle
of the $d+1$-dimensional real vector bundle
$\ss (\beta)$, where $\beta$ is the canonical
$2$-dimensional $\Kk$-vector bundle over $G_2^\Kk (\xi )$.
The bundle $\lambda$ is the real Hopf line bundle,
and $\alpha$ over $Q(\xi )$ is the orthogonal complement of
$\lambda$ in the lift of $\ss (\beta )$ to $P(\ss (\beta))$.

For any $\Kk$,
we now look only at $\Ff_2$-cohomology; the class $w_i^\Kk (\xi )$
reduces (mod $2$) to the Stiefel-Whitney class $w_{id}(\xi )$
(and $w_j(\xi )=0$ if $j$ is not divisible by $d$).
The description of the cohomology of the Grassmannian involves
polynomials $p_i(Y,Z)\in\Ff_2[Y,Z]$, $i=0,1,\ldots\,$, 
in indeterminates $Y$ of degree $d$ and $Z$ of degree $2d$,
defined by $p_0(Y,Z)=1$, $p_1(Y,Z)=Y$, and
$$
p_{i+1}(Y,Z)=Yp_i(Y,Z)+Zp_{i-1}(Y,Z) \text {\quad for $i\geq 1$.}
$$
\Lem{The $\Ff_2$-cohomology of $G_2^\Kk (\xi )$ is described as
$$
H^*(G_2^\Kk (\xi );\, \Ff_2)=H^*(B;\, \Ff_2)[Y,Z]/(p_{n}^\xi (Y,Z),\,
 Zp_{n-1}^\xi (Y,Z)+w_{(n+1)d} (\xi )),
$$
where 
$$
p_i^\xi (Y,Z)=p_i(Y,Z)+p_{i-1}(Y,Z)w_d (\xi )+\ldots
+p_1(Y,Z)w_{(i-1)d}(\xi )+w_{id} (\xi ),
$$
for $i\leq n$.
The Stiefel-Whitney classes of $\beta$ are
$Y=w_d (\beta)$, $Z=w_{2d}(\beta )$.
}
\begin{proof}
The relations come from
$$
(1+Y+Z)(1+w_d (\beta^\perp)+\ldots +w_{(n-1)d} (\beta^\perp))
=1+w_d(\xi )+\ldots +w_{(n+1)d} (\xi),
$$
where $\beta^\perp$ is the orthogonal complement of $\beta$ in 
the pullback of $\xi$.
So 
$$
1+w_d (\beta^\perp)+\ldots +w_{(n-1)d}(\beta^\perp)=
(1+\ldots +p_i(Y,Z)+\ldots )
(1+w_d (\xi )+\ldots +w_{(n+1)d}(\xi)),
$$
$$
p_{n}(Y,Z)+p_{n-1}(Y,Z)w_d(\xi ) +\ldots +
p_1(Y,Z)w_{(n-1)d} (\xi )
+w_{nd}(\xi )=0,
$$
and
$$
p_{n+1}(Y,Z)+p_{n}(Y,Z)w_d (\xi ) +\ldots +
p_1(Y,Z)w_{nd} (\xi )
+w_{(n+1)d}(\xi )=0,
$$
or
$$
Z(p_{n-1}(Y,Z)+p_{n-2}(Y,Z)w_d (\xi ) +\ldots +
w_{(n-1)d} (\xi ))
+w_{(n+1)d} (\xi )=0.
$$
Write $A=H^*(B;\, \Ff_2)$.
So we certainly have an $A$-homomorphism
$$
M= A[Y,Z]/(p_{n}^\xi (Y,Z),\,
 Zp_{n-1}^\xi (Y,Z)+w_{(n+1)d} (\xi )) \to H^*(G_2^\Kk (\xi );\, \Ff_2).
$$
By an application of the Leray-Hirsch Theorem, 
there is a finitely generated free
$A$-submodule $N\subseteq M$ which maps isomorphically
onto $H^*(G_2^\Kk (\xi );\, \Ff_2)$.
We are assuming that $B$ is connected, so that
$I=\tilde H^*(B;\, \Ff_2)$ is a nilpotent ideal
with $A/I= \Ff_2$.
By considering the restriction to a fibre, we
see that $M=IM+N$.
Hence, by Nakayama's lemma, we have $M=N$.
\end{proof}
\Prop{{\rm (Feder \cite{feder} for $B=*$).}
The $\Ff_2$-cohomology ring of $Q(\xi )$ is
$H^*(Q(\xi );\,\Ff_2)\,=$
$$
H^*(B;\,\Ff_2)[X,Y,Z]/
(X(X^d+Y), \,
p_{n}^\xi (Y,Z),\, Zp_{n-1}^\xi (Y,Z)+w_{(n+1)d} (\xi )),
$$
where the generators represent the $\Ff_2$-cohomology Euler classes
as: $X=e(\lambda )$, $Y=e(\alpha )+ e(\lambda )^d$, $Z=e(\beta )$.

In particular, $H^*(Q(\xi );\,\Ff_2)$
is free as a module over $H^*(G_2^\Kk (\xi );\,\Ff_2)$ with basis
$1,\, e(\lambda ),\,\ldots , e(\lambda)^d$.
}
\begin{proof}
For 
$$
H^*(P(\ss (\beta));\,\Ff_2)=H^*(G_2^\Kk (\xi );\,\Ff_2)[X]/(X^{d+1}+
w_1(\ss (\beta ) )X^d+\ldots +w_{d+1}(\ss (\beta ))),
$$
and $w_i(\ss(\beta ))=0$, if $0<i<d$,
$w_d(\ss (\beta ) )=w_d(\beta )$.
(We can compute when $\beta$ is a sum of two $\Kk$-line bundles.)

Since $\lambda\oplus\alpha \iso \ss (\beta)$, we have
$w_d(\alpha )=w_d(\ss (\beta ))+ w_1(\lambda )^d$.
\end{proof}
\Cor{\label{one}
If, for some $k\geq 1$, 
$e(\tilde\alpha )^{k-1}\in H^{(k-1)d}(\tilde Q(\xi );\,\Ff_2)$ 
is non-zero, then 
$e(\alpha )^{k}\in H^{kd}(Q(\xi );\,\Ff_2)$ is non-zero.
}
\begin{proof}
Since $(X^d+Y)^k -(X^d+Y)Y^{k-1}$ is divisible by $X(X^d+Y)$,
$e(\alpha )^k =w_d(\beta )^k\cdot 1+ w_d(\beta )^{k-1}\cdot e(\lambda)$
is non-zero if and only if $w_d(\beta )^{k-1}$ is non-zero.
But $w_d(\beta )=e(\alpha )+e(\lambda )^d$ lifts to $e(\tilde\alpha )$.
\end{proof}
\Ex{(Going back to the 1957 paper of Peterson \cite{fp}).
For $B=*$, $\xi =V$, $\Kk =\Rr$, $n=2^r$, 
condition (4) fails if $k=2^{r+1}-1=2n-1$,
but (1) holds if $k=2^{r+1}=2n$.
}
\begin{proof}
It is enough to show that the $\Ff_2$-cohomology Euler class
$e(\alpha )^{2^{r+1}-1}$ is non-zero. 
By Corollary \ref{one} this is true
if $e(\tilde\alpha )^{2^{r+1}-2}\not=0$, and this follows
from Example \ref{milnor}.
\end{proof}
\Sect{\label{rank}
Vector bundle homomorphisms of rank greater than $1$}
In this section we prove that,
if either (a) $n>1$ and $\dim B < (k-1)(n+1)$
or (b) $n=1$ and $\dim B <2(k-1)-2$, then
condition (4) in Theorem \ref{symm}
that $\gamma (\zeta )^k=0\in\omega^0(P(\xi );\, -k\zeta )$
implies condition (0) that there is a vector bundle map
$B\times\Rr^k\to\xi$ with rank greater than $1$ at each point.

\begin{proof}
We look first at the case (b): $n=1$ and $\dim B<2(k-1)-2$.
Existence of a map $B\times \Rr^k\to\xi$ with rank $\geq 2$ at each point, that is, a surjective map, is equivalent by duality to
existence of a bundle monomorphism $\xi\into B\times\Rr^k$.
In the stable range $\dim B <2(k-2)$ a monomorphism exists
if and only if 
$\gamma (\eta )^k\in\omega^0(P(\xi );\, -\eta\otimes\Rr^k)$ is zero.
But the involution of $P(\xi )$ taking a line to its orthogonal complement in the $2$-dimensional bundle $\xi$ 
maps $\eta$ to $\zeta$ and so $\gamma (\eta)$ to $\gamma (\zeta )$. 

For case (a) with $k=2$ and $\dim B <n+1$, that is, $\dim B\leq n$,
we can give a cohomological argument. We show that if
$w_n(\zeta )^2=0$, then $\xi$ admits a trivial summand
$B\times\Rr^2$. It suffices to prove that $w_n(\xi )=0$.
Now $H^*(P(\xi );\,\Ff_2)=H^*(B;\,\Ff_2)[t]/
(t^{n+1}+w_1(\xi )t^n+\ldots +w_n(\xi )t)$,
because $\dim B\leq n$. And $w_n(\zeta )
=t^n+w_1(\xi )t^{n-1}+\ldots +w_n(\xi )$.
Hence $tw_n(\zeta )=0$ and so $w_n(\zeta )^2=w_n(\xi )w_n(\zeta )
=w_n(\xi )t^n$. (For dimensional reaons,
$w_n(\xi )w_i(\xi )=0$ when $i\geq 1$.)
If $w_n(\zeta )^2=0$, it follows that $w_n(\xi )=0$.

Now consider the main case (a) with $k > 2$.
The argument which follows
can be understood as a special case of Koschorke's theory 
in \cite[Existence theorem 3.1]{ulrich}.

For fibrewise pointed spaces $X\to B$ and $Y\to B$ over $B$,
we write $\omega^0_B\{ X;\, Y\}$, as in \cite[II: Chapter 1]{ioan},
for the group of fibrewise stable maps from $X$ to $Y$. In general,
fibrewise constructions over $B$ are indicated by a subscript `$B$'
as: `$+B$' adding a disjoint basepoint in each fibre,
`$/_B$' forming the fibrewise topological cofibre (by collapsing a
subspace to a point), or the fibrewise Thom space of a vector bundle.

The sphere bundle $S(\Hom (\Rr^k,\eta ))$ over
$P(\xi )$ is included, fibrewise over $B$, 
as a fibrewise submanifold, $Z$ say, in
$S(\Hom (\Rr^k,\xi ))$ of dimension $n+k-1$.
Its fibre $Z_x$ at $x\in B$ is the closed manifold of
linear maps $\Rr^k\to\xi_x$ in $S(\Hom (\Rr^k,\xi_x))$
with rank equal to $1$.
The fibrewise normal bundle $\nu$ of $Z\subseteq
S(\Hom (\Rr^k,\xi ))$ has dimension $k(n+1)-1-(n+k-1)
=(k-1)n$.
(More precisely, $\nu$ is the cokernel of 
an inclusion of $\Hom (\eta ,\zeta )$ lifted to $Z$ into 
the pullback of $\Hom (\Rr^k,\zeta )$.)

Choose a fibrewise tubular neigbourhood $D(\nu )
\into S(\Hom (\Rr^k,\xi ))$ and let $W\to B$ be the (closed)
complement of the open tubular neighbourhood $D(\nu )-S(\nu)$.
We aim to show that $W\to B$ has a cross-section.
A section will give at each point $x\in B$ a linear map
$\Rr^k\to\xi_x$ with rank greater than $1$.

The stable homotopy exact sequence of the pair 
$(S(\Hom (\Rr^k,\xi ),W)$ over $B$ appears as the lefthand column
of the diagram:
$$
\begin{matrix}
\omega^0_B\{ B\times S^0;\, W_{+B}\} &\Rarr{c}{}&
\omega_B^0\{ B\times S^0;\, B\times S^0\}\\
\Darr{}{}&&\Darr{=}{}\\
\omega^0_B\{ B\times S^0;\, S(\Hom (\Rr^k,\xi ))_{+B}\} &
\Rarr{c}{}&\omega_B^0\{ B\times S^0;\, B\times S^0\}\\
\Darr{}{}&&\\
\omega^0_B\{ B\times S^0;\, S(\Hom (\Rr^k,\eta ))_B^\nu\} 
&\Rarr{\,\,\iso\,\,}{\phantom{\rm onto}}&
\omega^0_B\{ B\times S^0;\, S(\Hom (\Rr^k,\xi ))/_BW\} 
\end{matrix}
$$
The maps $c$ are induced by the projection of the fibres
of $W$ or $S(\Hom (\Rr^k,\xi ))$ to a point.
The (excision) isomorphism, involving the fibrewise Thom
space 
$$
D(\nu )/_B S(\nu )=S(\Hom (\Rr^k,\eta ))^\nu_B
$$ 
of $\nu$, is induced by
the inclusion $(D(\nu ),S(\nu ))\into (S(\Hom(\Rr^k,\xi ), W)$.

Let $\pi : P(\xi )\to B$ denote the projection.
We have a fibrewise inclusion 
$\iota :S(\Hom (\Rr^k,\eta ))\into S(\Hom (\Rr^k,\pi^*\xi ))$ 
over $P(\xi )$ 

By duality over $B$ -- see, for example,
\cite[II: Section 12]{ioan} -- we can express the relevant
fibrewise stable homotopy groups as stable cohomotopy groups:
$$
\omega^0_B\{ B\times S^0;\, S(k\xi )_{+B}\} =
\omega^{-1}(S(k\xi );\, -k\xi )
$$
and
$$
\omega^0_B\{ B\times S^0;\, S(k\eta )^\nu_B\}
=\omega^{-1}(S(k\eta );\, -k\xi ).
$$
These fit into a commutative diagram of Gysin sequences
$$
\begin{matrix}
\omega^{-1}(S(k\xi );\, -k\xi )&\Rarr{\rm onto}{}&
\omega^0(D(k\xi ),S(k\xi );\, -k\xi )
&=&\omega^0(B)\\
\Darr{\pi^*}{}&&\Darr{\pi^*}{}&&\Darr{\pi^*}{}\\
\omega^{-1}(S(k\pi^*\xi );\, -k\xi )&\Rarr{\rm onto}{}&
\omega^0(D(k\pi^*\xi ),S(k\pi^*\xi );\, -k\xi )
&=&\omega^0(P(\xi ))\\
\Darr{\iota^*}{}&&\Darr{\iota^*}{}&&\Darr{\iota^*}{}\\
\omega^{-1}(S(k\eta );\, -k\xi )&\Rarr{\iso}{\phantom{\rm onto}}&
\omega^0(D(k\eta ),S(k\eta );\, -k\xi )
&=&\omega^0(P(\xi );\, -k\zeta )
\end{matrix}
$$
The maps marked are, respectively,
surjective, because $\dim B <\dim (k\xi )$
so that $\gamma (k\xi )=0$,
and bijective, because $\dim P(\xi )+1=\dim B +n +1
< \dim (k\xi )=k(n+1)$, that is, $\dim B <(k-1)(n+1)$,
so that the groups $\omega^{-1}(P(\xi ); -k\xi )$ and
$\omega^{0}(P(\xi ); -k\xi )$ are zero.
In the righthand column,
$\iota^*(1)=\gamma (k\zeta )=\gamma (\zeta )^k$.

Now suppose that $\gamma (\zeta )^k=0$.
Then there is some class $x\in \omega^{-1}(S(k\xi );\, -k\xi )$
that maps to $1\in\omega^0(B)$ and by $\iota^*\pi^*$ to
$0\in\omega^{-1}(S(k\eta);\,-k\xi )$,
This translates back by duality into the existence of a class
$x\in\omega^0_B\{ B\times S^0;\, S(k\xi )_{+B}\}$ that
maps to $1\in\omega^0_B\{ B\times S^0;\, B\times S^0\}$
and to $0\in\omega^0_B\{ B\times S^0;\, S(k\xi )/_BW\}$.

We conclude from the stable cohomotopy exact sequence in the
first diagram
that $1\in\omega^0_B\{ B\times S^0;\, B\times S^0\}$ lifts to
$\omega^0_B\{ B\times S^0;\, W_{+B}\}$. 
This says that the bundle $W\to B$ admits a `stable section'.

But we are in the stable range
$\dim B < 2(\dim \nu -1) =2((k-1)n-1)$,
because $(k-1)(n+1)\leq 2((k-1)n-1)$,
that is, $(k-1)(n-1)\geq 2$, since we are assuming that
$k>2$ (and $n>1$). 
Hence $W\to B$ has a section, as required.
\end{proof}

\end{document}